\documentclass[letterpaper, 10 pt, conference]{ieeeconf}

\IEEEoverridecommandlockouts

\usepackage{cite}
\usepackage{amsmath,amssymb,amsfonts}
\usepackage{mathtools}
\usepackage{algorithmic}
\usepackage{graphicx}
\usepackage{subcaption}
\usepackage{caption}
\usepackage{times}
\usepackage{epsfig}
\usepackage{float}
\usepackage{todonotes}
\usepackage{textcomp}
\usepackage{xcolor}
\usepackage{authblk}
\usepackage{tikz}

\newtheorem{theorem}{Theorem}[section]
\newtheorem{lemma}[theorem]{Lemma}
\newtheorem{proposition}[theorem]{Proposition}

\newtheorem{remark}[theorem]{Remark}
\begin{document}
	
	\title{\centering Traffic Flow Reconstruction from Limited Collected Data}
	
	\author{Nail Baloul, Amaury Hayat, Thibault Liard and Pierre Lissy%
		\thanks{N. Baloul, A. Hayat, and P. Lissy are with CERMICS, Ecole des Ponts, Champs-sur-Marne, France. 
			{\tt\small nail.baloul@enpc.fr, amaury.hayat@enpc.fr, pierre.lissy@enpc.fr}}%
		\thanks{T. Liard is with XLIM, Université de Limoges, Limoges, France. 
			{\tt\small thibault.liard@unilim.fr}}%
		
		\thanks{NB and PL are supported by ANR-22-CPJ2-0138-01 TL is supported by ANR-24-CE46-7295, and AH was supported by ANR-24-ERC S-010. This work was granted access to the HPC resources of IDRIS under the allocation 2023-AD011014527 made by GENCI.}
	}	\maketitle
	
	\begin{abstract}
		We propose an efficient method for reconstructing traffic density with low penetration rate of probe vehicles. Specifically, we rely on measuring only the initial and final positions of a small number of cars which are generated using microscopic dynamical systems. We then implement a machine learning algorithm from scratch to reconstruct the approximate traffic density. This approach leverages learning techniques to improve the accuracy of density reconstruction despite constraints in available data. For the sake of consistency, we will prove that, if only using data from dynamical systems, the approximate density predicted by our learned-based model converges to a well-known macroscopic traffic flow model when the number of vehicles approaches infinity.
		
	\end{abstract}
	
	
	\section{Introduction}
	
	Vehicular traffic is often described at two distinct levels, microscopic and macroscopic. The microscopic approach focuses on individual vehicles and their interactions. It considers the behavior of each driver capturing detailed dynamics such as velocity and acceleration. It enables to account for the impact of individual decisions on overall traffic flow. The macroscopic approach relies on a continuity assumption, describing traffic flow in terms of density, flow rate and mean speed. It treats traffic as a continuous fluid, governed by conservation principles such as the conservation of vehicles. It allows for large-scale traffic management. We refer to the survey papers \cite{piccoli_vehicular_2009,bellomo_modeling_2011,albi_vehicular_2019} for a general discussion about traffic models at different scales. Because the state of continuous systems or high-dimensional microscopic systems is impossible to fully monitor, being able to reconstruct it from limited traffic data is a paramount challenge that has attracted a strong attention in the last two decades \cite{herrera_evaluation_2010,seo_traffic_2017,barreau_learning-based_2021,barreau_physics-informed_2021}.
	
	We aim
	to construct a model that closely approximates the original (or intrinsic) traffic model while observing considerably fewer vehicles. This approach is motivated by the fact that monitoring all vehicles on a roadway is very challenging and often impractical due to technical limitations or resource constraints. 
	Our main contribution is to implement a machine learning algorithm that reconstructs traffic flow using only the initial and final positions of probe vehicles.  A key advantage of our approach is that it implicitly satisfies conservation laws without having to enforce them with additional constraints.
	This is achieved by using artificial data created by adjusting well-known traffic flow models and using activation functions in neural networks tailored to underlying physical problems.
	
	The paper is organized as follows.
	In section 2, we introduce the traffic models used throughout our study.
	Section 3 paves the way for the formulation of a constrained optimization problem that will be the core of our machine learning algorithm in section 4.
	Section 5 provides a theoretical guarantee of our model.
	In section 6 we will illustrate the efficiency of our learning-based method by numerical experiments.
	Finally, section 7 concludes the paper and suggests directions for future research.
	
	\section{Traffic flow models}
	
	\subsection{Microscopic model}
	We consider a system of $N+1$ vehicles moving on a one dimensional single-lane road. The first vehicle is designated as the leader while the remaining vehicles are referred to as followers. The position of each vehicle $i$ at time $t$ is denoted by $x_i^N(t)$.
	Each vehicle, regarded as a particle, moves at a speed 
	computed via a velocity map 
	$v(\cdot)$ and a distance to the vehicle immediately ahead of it. We assume that the leader has no vehicle in front of it and therefore maintains constant and maximal velocity denoted $v_{\max}$. Additionally, the cumulative length of the vehicles is given by $L$. In other words, $L=lN$, where $l$ is the average length of a single car. Then the Follow-the-Leader (FtL) microscopic 
	model is 
	\begin{equation}
		\begin{cases}
			\dot{x}^N_{N}(t)=v_{\max},&\quad t>0,
			\\
			\dot{x}^N_{i}(t)=v\left(\frac{L}{N\left(x^N_{i+1}(t)-x^N_{i}(t)\right)}\right),&\quad t>0,\\
			x^N_i(0)=\bar{x}_i,&\quad i=0,\dots, N,
		\end{cases}\label{2b}
	\end{equation}
	where $\bar{x}_i$ denotes the initial position of vehicle $i$. In \cite{francesco_rigorous_2017}, the authors established well-posedness of system \eqref{2b} under the condition that the initial positions satisfy $\bar{x}_0<\bar{x}_1<\dots<\bar{x}_{N-1}<\bar{x}_{N}$. This nonlinear model not only becomes high-dimensional and computationally expensive when dealing with a large number of vehicles but also assumes that drivers react solely to the vehicle immediately ahead, disregarding the influence of other vehicles and overall traffic conditions. Thus when we have a large number of vehicles it is often more adapted -at least computationally- to consider a macroscopic model as described in the following.
	
	\subsection{Macroscopic model}
	The Lighthill-Whitham-Richards (LWR) model describes traffic flow using the one-dimensional conservation law. Indeed, vehicles are treated as a continuous medium similar to particles in fluid which leads us to the initial value problem
	\begin{equation}
		\begin{cases}
			\frac{\partial \rho}{\partial t}(t, x) + \frac{\partial f(\rho)}{\partial x}(t, x) = 0, &\quad x\in\mathbb{R}, \quad t>0,\\
			\rho(0, x)=\bar{\rho}(x), &\quad x\in\mathbb{R},
		\end{cases}\label{1}
	\end{equation}
	where $\rho(t, x)$ is the traffic density at position $x$ and at time $t$, $f(\rho)$ is the flux function representing the flow rate of vehicles and $\bar{\rho}$ is the initial vehicular density. Typically, the flux function is expressed as $f(\rho) := \rho v(\rho)$, where $v(\rho)$ is representing the average velocity of vehicles at a given density. Obviously, this model does not capture the behavior of individual vehicles as it focuses on aggregate traffic flow variables such as density and flow rate. 
	
	The convergence analysis of microscopic FtL model \eqref{2b} towards macroscopic LWR model \eqref{1} has been extensively studied in recent years \cite{colombo_well_2007}, \cite{holden_continuum_2019}. The procedure involves considering an initial density $\bar{\rho}$ which is discretized to determine the initial positions of $N+1$ vehicles. We then allow FtL system \eqref{2b} to evolve according to its dynamics and compare the resulting solution with the one that satisfies LWR model \eqref{1} as the number of vehicles approaches infinity and the length of each car tends to zero.
	
	\subsection{Data-driven based model}
	
	In our scenario, it is crucial to note that, unlike previous works that often rely
	on a given initial density, we do not have access to this critical information. It represents 
	a significant challenge, as density is a key component for accurately predicting traffic flow. 
	
	In \cite{barreau_physics-informed_2021,barreau_learning-based_2021}, the authors considered probe vehicles regarded as mobile sensors within the traffic flow. These vehicles are governed by an Ordinary Differential Equation (ODE) system analogous to \eqref{2b} where traffic density evolves according to a Partial Differential Equation (PDE) related to \eqref{1}. To reconstruct density, vehicle positions $x_i(t)$, local density $\rho_i(t)=\rho(t, x_i(t))$ at vehicle locations and instantaneous speed $v_i(t)=\dot{x}_i(t)$ were needed in real time. This physics-informed learning method relies on explicit enforcement of conservation laws through PDE-derived Lagrangian terms in their optimization. Its main drawback is the lack of theoretical guarantee to prove that the reconstructed density converges to the conservation law.
	In \cite{inzunza_pinn_2023}, the authors adopted a Physical Informed Neural Network (PINN) strategy using density and flow measurements from synthetic data computed at fixed spatial locations. They illustrated that accuracy of their traffic state estimation improves when the number of locations increases.
	In \cite{shi_physics_2021}, the authors designed a Physics-Informed Deep Learning (PIDL) method that integrates machine learning to learn the Fundamental Diagram (FD) mapping from traffic density to flow or velocity.
	In \cite {herrera_evaluation_2010}, the authors conducted the Mobile Century field experiment where they 
	used GPS-enabled mobile phones as a cost-effective traffic monitoring system. The system collected continuous trajectory data, with 
	updates recorded every three seconds.
	Their results suggested that a $2-3\%$ penetration of GPS-equipped phones in the driver population could provide accurate measurements of traffic flow velocity. 
	
	In contrast, our approach requires fewer information than \cite{barreau_physics-informed_2021,barreau_learning-based_2021,inzunza_pinn_2023, herrera_evaluation_2010, shi_physics_2021} as it only assumes the knowledge of initial and final positions of 
	probe vehicles, 
	simplifying 
	data collection 
	and reducing the need for real-time measurements. 
	
	\section{Mathematical formulation}
	To formalize our model, we first adapt 
	\eqref{2b} to a 
	framework involving 
	an optimization variable and a finite time horizon.
	\subsection{Parametrized ODE system} 
	We consider the ODE system with finite time $T$
	\begin{equation}
		\begin{cases}
			\dot{x}^N_{n}(t)=v_{\max},&\quad t\in (0,T] ,
			\\
			\dot{x}^N_{i}(t)=v\left(\frac{\alpha^N_{i}L}{N\left(x^N_{i+1}(t)-x^N_{i}(t)\right)}\right),&\quad t\in (0,T], 	\\
			x^N_i(0)=\bar{x}_i,&\quad i=0,\dots, n,
		\end{cases}\label{3b}
	\end{equation}
	where $\alpha^N=(\alpha^N_0,\alpha^N_1,\dots,\alpha^N_{n-1})^\top\in\mathbb{R}^n$ is a parameter and $L$ is still representing the total length covered by all 
	vehicles in a bumper-to-bumper configuration. 
	In 
	system \eqref{3b}, 
	considered vehicles represent only a subset of the total involved. For $i=0,\dots,n-1$, $\bar{x}_i$ and $\bar{y}_i$ denote respectively the given initial and final positions of car $i$. The goal is to gather observations from a significantly reduced number $n$ of probe vehicles compared to the total number $N$. 
	
	In system \eqref{3b}, $\alpha_{i}^{N}$ stands for the number of vehicles in the segment $[x_i(\cdot), x_{i+1}(\cdot))$ between consecutive probe vehicles $i$ and $i+1$. Therefore, we impose physical conditions 
	\begin{equation}
		\alpha^N_i\in\left[1, \bar{z}_i\right],\quad i=0,\dots, n-1, 
		\quad\sum_{j=0}^{n-1}\alpha^N_j=N,
	\label{alpha_constraints}
\end{equation}
where 
	$\bar{z}_i\coloneqq\min\left\{N\left(\bar{x}_{i+1}-\bar{x}_{i}\right)/L, N\left(\bar{y}_{i+1}-\bar{y}_{i}\right)/L\right\}.$
The lower bound of 1 is trivial,
as the segment must contain at least one vehicle. The upper bound ensures that 
there is enough space between consecutive probe vehicles for this number of vehicles. Since $L/N$ is the length of a car, there cannot be more than $k$ cars on a segment, as 
$kL/N$ is the minimal possible space required by $k$ cars on the road.
Moreover, a key assumption of our model is that overtaking between probe vehicles is not allowed. In other words, similarly to  FtL system \eqref{2b}, the number of vehicles between consecutive probe vehicles remains constant over time which translates to $\alpha_{i}^{N}$ being constant with time.
To handle constraints \eqref{alpha_constraints} efficiently, we introduce the notation
\begin{equation}
	\mathcal{A}_N\coloneqq\left\{\alpha \in \mathbb{R}^n: \alpha^N_i\in\left[1, \bar{z}_i \right],  \sum_{j=0}^{n-1}\alpha^N_j=N\right\}\label{alpha_set}.
\end{equation}

The inclusion of $\alpha^N$ in ODE system \eqref{3b} preserves well-posedness of the system while allowing for a parametric adjustment to the vehicle dynamics. 
It will also play a role in establishing convergence to macroscopic model \eqref{1} which facilitates the transition from discrete, vehicle-level dynamics to continuous, density-level dynamics.\\
We choose from now on to drop the superscript when referring to trajectories for simplicity of notation.

\subsection{Global existence and uniqueness}

We assume that the velocity 
satisfies 
asumptions
\newcounter{vctr}
\renewcommand{\thevctr}{v\arabic{vctr}}

\begin{list}{(\thevctr)}{
		\usecounter{vctr}
		\setlength{\leftmargin}{2em}
		\setlength{\labelwidth}{2em}
		\setlength{\itemsep}{0.5ex}
	}
	
	\item \label{v1} $v \in C^1([0,+\infty))$,
	\item \label{v2} $v$ is decreasing on $[0,+\infty)$,
	\item \label{v3} $v(0) = v_{\max}$ for some $v_{\max} \in \mathbb{R}$.
\end{list}
Since parameter $\alpha^N$ is bounded and $v$ is smooth, it is straightforward to show local existence and uniqueness of the solution to ODE system \eqref{3b} from Picard-Lindel\"{o}f Theorem. To ensure that the solution exists globally, we have to prove that the distance between two consecutive probe vehicles never vanish. The result has been shown by Di Francesco and Rosini in \cite{francesco_rigorous_2017} for the classical FtL model \eqref{1} involving all vehicles in the system. We briefly outline it in our context, where the same principle can be proven to hold.

\begin{lemma}
	\label{lemma1}
	Let $\left(x_0(\cdot), \dots, x_n(\cdot)\right)$ be the solution of ODE system \eqref{3b} and $v$ satisfy hypotheses (\ref{v1})-(\ref{v3}). Then the discrete maximum principle holds; for all $i=0,\dots, n-1$ and for all $t\in [0,T]$,
	\begin{equation}
		\frac{\alpha^N_iL}{NM}\leq x_{i+1}(t)- x_i(t)\leq \bar{x}_n-\bar{x}_0+\left(v_{\max}-v(M)\right)t,
	\end{equation}
	where $M\coloneqq   \underset{i=0,\dots, n-1}{\max}\left(\frac{\alpha^N_i L}{N(\bar{x}_{i+1}-\bar{x}_{i})}\right)$ denotes the maximum discrete density at initial time $t=0$.
\end{lemma}

\subsection{Discrete density}
To clearly articulate the relationship between the state of ODE system \eqref{3b} and the unknown density in equation \eqref{1} we introduce for $i=0,\dots, n-1$ the discrete density as
\begin{equation}
	\rho_i^N(t) \coloneqq  \frac{\alpha^N_iL}{N(x_{i+1}(t)-x_{i}(t))},\quad t\in [0,T],\label{4}
\end{equation}
where the trajectories $x(\cdot)$ of the probe vehicles are solutions to system \eqref{3b}. $\rho_i^N(\cdot)$ represents the local density in the interval between consecutive probe vehicles localised at $x_i(\cdot)$ and $x_{i+1}(\cdot)$.
Based on \eqref{4}, we define for $x\in\mathbb{R}$ the piecewise constant Eulerian discrete density as
\begin{equation}
	\rho^N(t,x)\coloneqq   \sum_{i=0}^{N-1}\rho^N_i(t)\chi_{[x_i(t),x_{i+1}(t))}(x),\quad t\in[0,T],\label{9a}
\end{equation}
where $\chi_A$ is the characteristic function of the set $A$.
As proved in \cite{holden_continuum_2019}, $\rho^N$ in \eqref{9a} can be seen as a discrete approximation of the solution to initial value problem \eqref{1}.

\section{Learning-based density estimation}

Our main focus is to reconstruct the traffic density only using initial and final positions of probe vehicles. To that end, we will solve an optimization problem that will allow us to construct a traffic discrete density $\rho^N$ in \eqref{9a} that converges to the solution of LWR model \eqref{1} when $N$ tends to infinity.

\subsection{ODE constrained optimization problem}

We address our optimization problem with an innovation which contrasts with physics-informed methods. It lies in deliberately avoiding direct incorporation of the PDE into the system due to the lack of theoretical guarantee induced by these methods. 
We instead leverage the established convergence under certain conditions of FtL microscopic model \eqref{2b} to LWR macroscopic model \eqref{1}. This convergence implicitly ensures mass conservation without requiring explicit PDE constraints, thus maintaining flexibility in the reconstructed solution space while preserving physical consistency. 
Enforcing directly the PDE constraints  would not pose a computational challenge but would instead create a data dependency. Indeed, it would require time-varying density measurements at specific points in space, which we exclude from our problem set-up since our methodology assumes no density measurements.
We first 
compute the leader's trajectory 
\begin{equation*}
	x_n(t)=v_{\max}t+\bar{x}_n,\quad t\in [0,T].\label{13}
\end{equation*}
We then introduce for $\alpha^N \in\mathcal{A}_N$ the matrix $W_{\alpha^N}$ that accounts for interaction between probe vehicles
\begin{equation}
	\begin{cases}
		(W_{\alpha^N})_{i,i}\coloneqq  -\frac{N}{\alpha^N_iL}, &i=0,\dots,n-1,\\
		(W_{\alpha^N})_{i,i+1}\coloneqq  \frac{N}{\alpha^N_iL}, &i=1,\dots,n-2,\\
		(W_{\alpha^N})_{i,j}\coloneqq  0, &\text{otherwise} ,
	\end{cases}\label{15}
\end{equation} 
and function $b_{\alpha^N}(\cdot)$ defined as
\begin{equation}
	b_{\alpha^N}(t)=\left(0,\dots,0,b_{n-1}(t)\right)^\top\in\mathbb{R}^{n},\label{16}
\end{equation}
where $b_{n-1}(t)\coloneqq  N / \alpha^N_{n-1}L\left(v_{\max}t+\bar{x}_{n}\right)$ accounts for the influence of the leader towards its followers.
Setting $V$ as
\begin{equation}
	V(z) \coloneqq   {v\left(z^{-1}\right)},\quad z\in\mathbb{R},\label{8}
\end{equation} 
we can rewrite system \eqref{3b} as
\begin{equation}
	\begin{cases}
		\dot{x}(t)=V\left(W_{\alpha^N}x(t)+b_{\alpha^N}(t)\right),&\quad t\in (0,T], \\x(0)=(\bar{x}_0,\dots, \bar{x}_{n-1})^\top,
	\end{cases}\label{17}
\end{equation}
where the state variable $x(\cdot)\coloneqq  (x_0(\cdot),\dots,x_{n-1}(\cdot))^\top$ represents the trajectories of the probe followers.


We emphasize that initial density $\bar{\rho}$ of the entire fleet is unknown, as we do not have access to the starting positions of all vehicles. Consequently, the discretized initial density $\bar{\rho}^N$ which would typically depend on initial positions of all vehicles, cannot be determined. Equivalently, the ground truth value of $\alpha^N$ is unknown as it relies on information that is not available to us. This limitation highlights the challenge of reconstructing the system's behavior based on our limited data collection.
From  Cauchy problem \eqref{17} and feasible set \eqref{alpha_set} we write our ODE constrained optimization problem as 
\begin{subequations}
	\begin{align}
		&\underset{\alpha^N}{\text{minimize}} \quad \frac{1}{2}\left\lVert x(T) - \bar{y}\right\rVert^2\label{eq:objective}\\
		\text{s.t.} \quad 
		& \dot{x}(t)=V\left(W_{\alpha^N}x(t)+b_{\alpha^N}(t)\right),\quad t\in[0,T], \label{eq:vehicle_dynamics}\\
		& x(0)=\bar{x},\label{eq:initial_conditions} \\
		& \alpha^N \in \mathcal{A}_N.\label{eq:alpha_constraints}
	\end{align}\label{optim}
\end{subequations}
The objective function \eqref{eq:objective} 
is a final position matching term which minimizes the discrepancy between the observed final positions $\bar{y}$ of probe vehicles and the corresponding final positions $x(T)$ predicted by the microscopic ODE model.  Moreover, the problem involves constraints related to dynamics of the probe vehicles \eqref{eq:vehicle_dynamics}-\eqref{eq:initial_conditions} whose velocities depend on inter-vehicle spacing influenced by $\alpha^N$ that must satisfy \eqref{eq:alpha_constraints}.
This formulation enables traffic density reconstruction from sparse observations by enforcing realistic vehicle dynamics, without directly incorporating PDE constraints. 

Machine learning is effective at finding complex data patterns.
Our method is based on the use of a neural network designed to understand the dynamics of traffic by breaking down the process into small time steps.

\subsection{Dataset generation}
Dataset consist of artificial data based on simulated FtL dynamics.
Precisely, we let evolve $N+1$ vehicles using classical FtL model \eqref{2b} until time $T$ is reached.
$10\%$ of the total fleet serve as probe vehicles for training data.
These vehicles provide initial and final position data used in the machine learning training process.
The probe vehicles represent a balanced sample of the overall traffic dynamics.
An additional $2.5\%$ of the total fleet are selected for testing purposes so that $80\%$ of the artificial data are for training and the remaining $20\%$ are test data.
Test data are strictly reserved for testing the model and are not used in any part of the training process. It is crucial to emphasize that the optimization network never ``sees" or interacts with the initial and final positions of test vehicles. This separation is fundamental to our methodology, as it allows for an unbiased assessment of the model's performance.
This choice enables us to evaluate whether the model has truly learned the underlying physics of traffic flow, rather than merely memorizing training data.


\subsection{Learning architecture}
We have designed the network architecture using a residual network (ResNet) approach. Each residual block in our ResNet corresponds to a single time step in the simulation. 
The network starts with an input representing the initial positions of probe vehicles.
This state is then propagated through time using a single, repeating residual block. The number of times this block is applied corresponds to the number of time steps in the simulation. 
Starting from initial traffic positions, the network predicts the next traffic state based on the current one. This repeated process allows step-by-step simulation of traffic flow until final time is reached. The structure described mirrors the first-order Euler discretization for ODE system \eqref{17}. 
The residual block incorporates physics-based principles of traffic flow. The weights \eqref{15} and biases \eqref{16} of the network are not independent variables, but are instead functions of $\alpha$. Nonlinear map \eqref{8} acts as a physics-grounded activation function. Optimizing $\alpha$ determine the entire network's behavior.

\usetikzlibrary{positioning,fit}
\definecolor{forward}{RGB}{0,102,204}
\definecolor{backward}{RGB}{204,0,0}

\begin{figure}[!t]
	\centering
	\begin{tikzpicture}[
		node distance=0.8cm and 0.8cm,
		font=\footnotesize, 
		forward/.style={color=forward, thick},
		backward/.style={color=backward, thick, dashed}
		]
		\node (dynamics) [draw=forward, rounded corners, align=center, inner sep=2pt] {NN Parameter $\alpha$};
		\node (fixed step) [draw=forward, rounded corners, below=of dynamics, yshift=0.2cm, inner sep=2pt] {Discretization step $\Delta t$};
		\node (residual) [draw=forward, rounded corners, below=of fixed step, yshift=0.2cm, inner sep=2pt] {Dynamics $f_\alpha(x_k, i\Delta t)$};
		\node (state) [draw=forward, rounded corners, below=of residual, yshift=0.2cm, inner sep=2pt] {State $x_{k+1} = x_k + f_\alpha(x_t, i\Delta t) \Delta t$};
		\node (repeat) [draw=forward, rounded corners, below=of state, yshift=0.4cm, inner sep=2pt] {Repeat for all $k = i\Delta t$};
		\node (final state) [draw=forward, rounded corners, below=of state, yshift=-0.3cm, inner sep=2pt] {Final state $x(T)$};
		\node (loss) [draw=backward, rounded corners, right=of final state, xshift=1.4cm, inner sep=2pt] {Loss $\mathcal{L}^{\mathrm{train}}(x(T), \bar{y})$};
		\node (optim) [draw=backward, rounded corners, above=of loss, inner sep=2pt] {$\alpha \leftarrow \alpha - \eta\nabla_\alpha \mathcal{L}$};
		
		\node[fit=(residual)(state), draw, dashed, inner sep=0.3cm] {};
		
		\draw[->, forward] (dynamics) -- (fixed step);
		\draw[->, forward] (fixed step) -- (residual) node[midway, left] {\scriptsize Residual Block $i$};
		\draw[->, forward] (residual) -- (state);
		\draw[->, forward] (repeat) -- (final state);
		\draw[->, forward] (final state.east) -- (loss.west) node[midway, above, font=\scriptsize] {Prediction};
		
		\draw[->, backward] (loss) -- (optim) node[midway, right, font=\scriptsize] {Auto-diff.};
		\draw[->, backward] (optim.north) |- (dynamics.east) node[pos=0.70, above, font=\scriptsize] {Parameter update};
		
		\node[anchor=north west, font=\scriptsize] at ([xshift=-0.2cm,yshift=-0.4cm]current bounding box.south west) {
			\begin{tabular}{@{}l l}
				\tikz\draw[->, forward] (0,0) -- (0.4,0); & Forward process \\
				\tikz\draw[->, backward] (0,0) -- (0.4,0); & Backward propagation
			\end{tabular}
		};
		
		\draw[black] (current bounding box.south west) rectangle (current bounding box.north east);
		
	\end{tikzpicture}
	\caption{Learning procedure}
	\label{fig:learning-procedure}
	\vspace{-\baselineskip}
\end{figure}
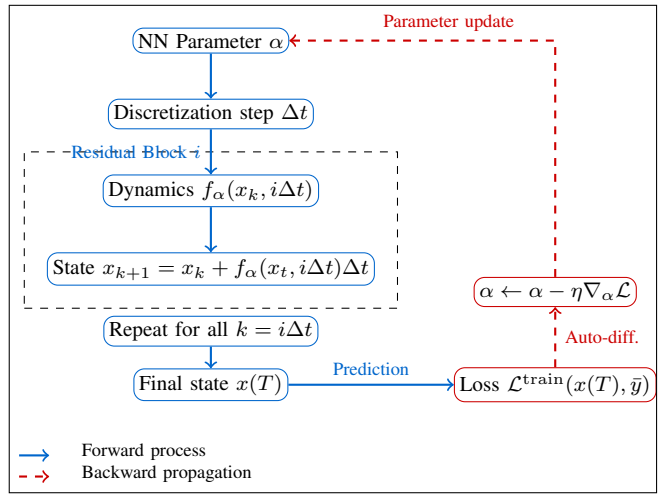

\subsection{Learning procedure}

Fig.~\ref{fig:learning-procedure} illustrates the overall learning procedure. 
Our task is to minimize the cost between predicted and observed final positions.
Predictions made by the neural network follow traffic rules defined by dynamics \eqref{eq:vehicle_dynamics} parameterized by $\alpha$. Indeed, the final traffic state $x^{\alpha}(T) $ reconstructed through our ResNet architecture is compared to the observed training data $\bar{y}$.
The training loss derived from \eqref{eq:objective} writes
$\mathcal{L}^{\mathrm{train}}=\frac{1}{n}\sum_{j=0}^{n}\lvert x_j^{\alpha}(T)-\bar{y}_j\rvert^2$ 
and computes the error between this final predicted state and the final observed state. Optimization parameter $\alpha$ is adjusted using backpropagation to minimize prediction errors. 
This iterative process of prediction, comparison and adjustment enables the model to capture the complex dynamics of traffic flow.

After training, we recover from $\alpha$ the discrete vehicle densities \eqref{4}. From these, we construct piecewise constant Euler discrete density $\rho^N$ in the form of \eqref{9a}  and use it to simulate on our test data.  In particular, we solve system \begin{equation}
\begin{cases}
	\dot{x}_{i}(t)=v\left(\rho^N(t,x_i(t)^{_+})\right),&\quad t\in (0,T], \\
	x_i(0)=\bar{x}_i^N&\quad i=0,\dots, n_{\mathrm{test}},
\end{cases}\label{ode_test}
\end{equation}
where $\rho\left(t,x_i(t)^+\right)\coloneqq 
\lim\limits_{x\rightarrow x_i(t);x>x_i(t)}
\rho\left(t,x\right)$ which indicates that vehicle $i$ is affected by the density located in front and where the density $\rho^N$ is the approximate density obtained from training. We let evolve system \eqref{ode_test} from initial test positions $\bar{x}_i^{n_{\mathrm{test}}}$ and measure the testing error $\mathcal{L}^{\mathrm{test}}=\frac{1}{n_{\mathrm{test}}}\sum_{j=0}^{n_{\mathrm{test}}}\lvert x_j^{\alpha}(T)-\bar{y}_j\rvert^2$ between simulated final positions and observed test data $\bar{y}_i^{n_{\mathrm{test}}}$. 

\section{Convergence of the model}

In this section we prove that, if only using data from dynamical systems, the approximate density $\rho^N$ constructed by our machine learning model converges to the solution of LWR 
\eqref{1} when the number of vehicles approaches infinity.

Although our setting incorporates probe vehicles and parameter $\alpha$, we can adapt most of the results established 
in \cite{francesco_rigorous_2017}. We recall that the authors showed the convergence of the solution to microscopic model \eqref{2b} to the unique entropic solution to macroscopic model \eqref{1} when $N\rightarrow+\infty$. 
In our context, the main challenge lies in imposing a condition on the distribution of $\alpha$ which would guarantee convergence. 

We demonstrate that the discrete initial density $\rho^N(0,\cdot)$ converges to the initial condition  $\bar{\rho}$ in LWR model \eqref{1} under this additional assumption. 
We make the important observation that by construction the initial traffic density must satisfy for $i=0,\dots, n-1$
\begin{equation}
\bar{x}_{i+1}=\sup\left\{x\in \mathbb{R}:\int_{\bar{x}_i}^{x}\bar{\rho}(y) dy\leq \frac{\alpha_{i} L}{N}\right\}.\label{atom}
\end{equation} 
Indeed, although we do not have access to the ground-truth initial car density $\bar{\rho}$,  we do know that initial positions $\bar{x}_i$ of our probe vehicles verify \eqref{atom} which states that the number of unobserved vehicles between $[x_i, x_ {i+1})$ is given by $\alpha_i$.

The following result inspired from \cite[Proposition 4]{francesco_rigorous_2017} ensures that the discretization aligns consistently with the true initial density when $N$ tends to infinity. 
To simplify notations, let for $t\in[0,T]$, $\rho(t)\coloneqq \rho(t,\cdot)$ .
We require the Wasserstein distance defined in \cite{francesco_rigorous_2017} by
\begin{equation}
\mathrm{W}_{L,1}(f, g) = \|f((-\infty,\cdot])-g((-\infty,\cdot)\|_{L^{1}(\mathbb{R},\mathbb{R})}.\label{wasserstein}
\end{equation}
We refer readers to \cite[Section 2.3]{francesco_rigorous_2017} for a rigorous introduction to this concept. Then, one can prove the following.

\begin{proposition} Let $\bar{\rho}$ satisfies \eqref{atom} and assume that \ \begin{equation}
	\max_{i=0,\dots, n-1}\alpha^N_i = o(N).\label{alpha_assumpt}
\end{equation}
Then, 
the sequence $(\rho^N(0))_{n \in \mathbb{N}}$
converges 
to $\bar{\rho}$ in the sense of the $\mathrm{W}_{L, 1}$
-Wasserstein distance in \eqref{wasserstein}.
\end{proposition}
\begin{remark}  
A particular case of assumption \eqref{alpha_assumpt} is when $\max_{i=0,\dots, n-1}\alpha^N_i \leq  CN\log(N)^{-1}$ for some $C>0$. 
This common choice
ensures controlled growth of $\alpha_N$.
\end{remark}
Moreover, by leveraging the expression of discrete density \eqref{9a}, 
we can generalize the convergence to the entropy solution of conservation law \eqref{1}, referring to the methodologies from abovementioned work. Specifically, \cite[Theorem 3]{francesco_rigorous_2017} which asserts the convergence of $\rho^N$ to the unique entropy solution of \eqref{1} remains applicable in our case, requiring only minor modifications to the original arguments.

\section{Numerical experiments}

We present some results of our numerical simulations using the abovementioned training procedure. 
We consider two distinct traffic scenarios where the maximum allowable speed of traffic is $v_{\max} = 120$ km/h and the maximum traffic density is $\rho_{\max} = 200$ cars/km. We consider a special instance of function $v$, known as the Greenshields velocity 
$v(\rho)=v_{\max}\max\left\{1-\rho / \rho_{max},0\right\},\quad \rho\in [0,\rho_{max} ].$ 
A time period of $0.1$ corresponds to a real-world duration of $6$ minutes. In the following illustrations, the origin represents the initial position $\bar{x}_0 = 0$ of the last follower at $t = 0$. The right end of the x-axis corresponds to the leader's position at the final time $T$. The initial domain refers to the road portion occupied by all probe vehicles at the simulation start.
The first scenario captures a 
traffic pattern characterized by alternating regions of congestion and free flow, represented by multiple waves.

Let $\rho(0,x)=0.6\rho_{\max} +0.3\rho_{\max}\sin(2\pi k x)$ be the initial density, where $k>0$ is the number of waves. We set $k = 3$.


\begin{figure}
\begin{subfigure}[b]{0.23\textwidth}
	\centering
	\includegraphics[width=\textwidth]{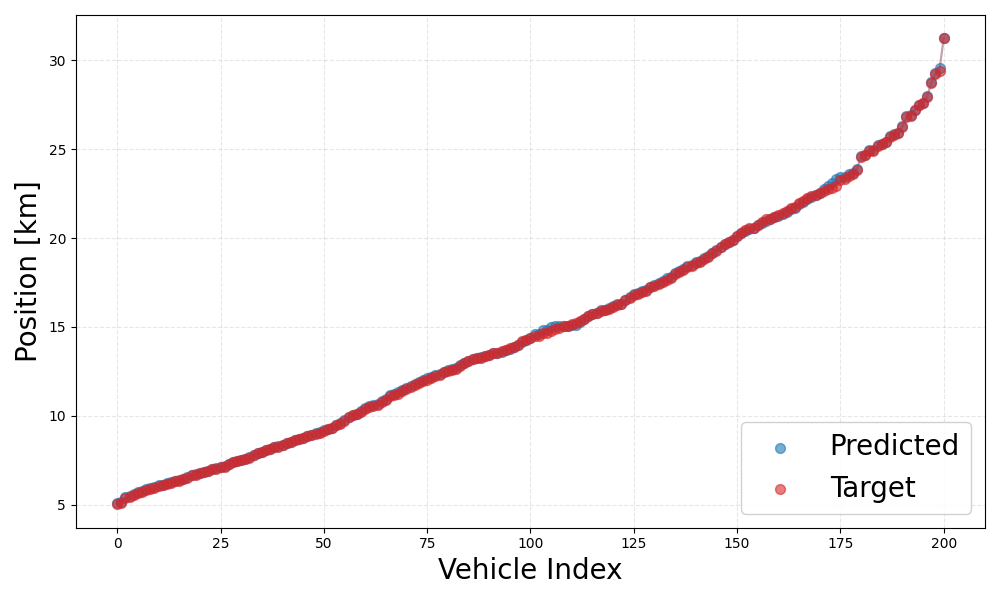}
	\caption{NN learning using 5,000 epochs}
	\label{fig:predicted_vs_target_36}
\end{subfigure}
\hfill
\begin{subfigure}[b]{0.23\textwidth}
	\centering
	\includegraphics[width=\textwidth]{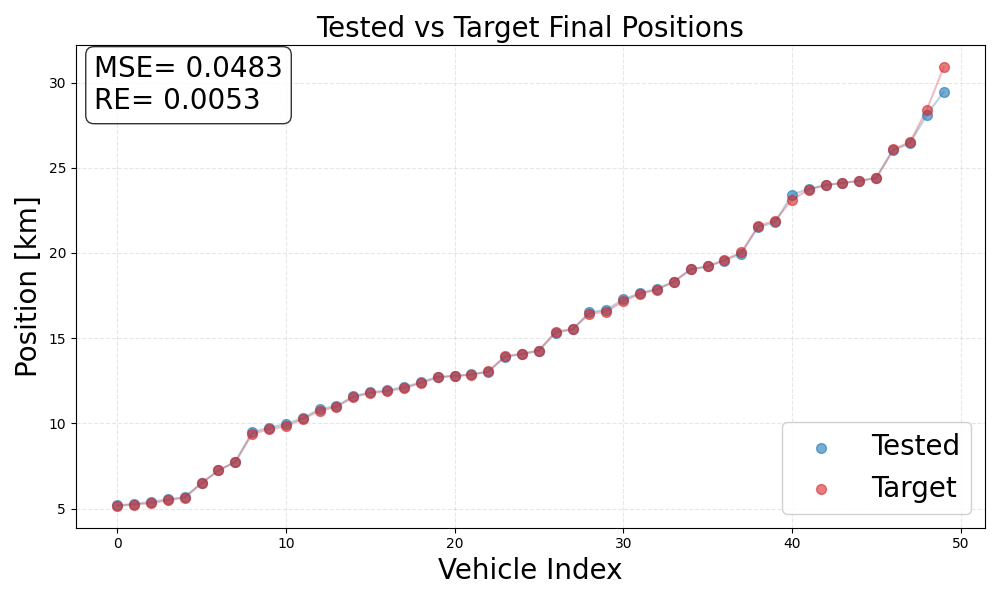}
	\caption{Testing using predicted density}
	\label{fig:tested_vs_target_36}
\end{subfigure}
\caption{Comparison between predictions and data with $N=2000$}
\label{fig:stop_and_go_recon_36}
\end{figure}

\begin{figure}
\begin{subfigure}[b]{0.23\textwidth}
	\centering
	\includegraphics[width=\textwidth]{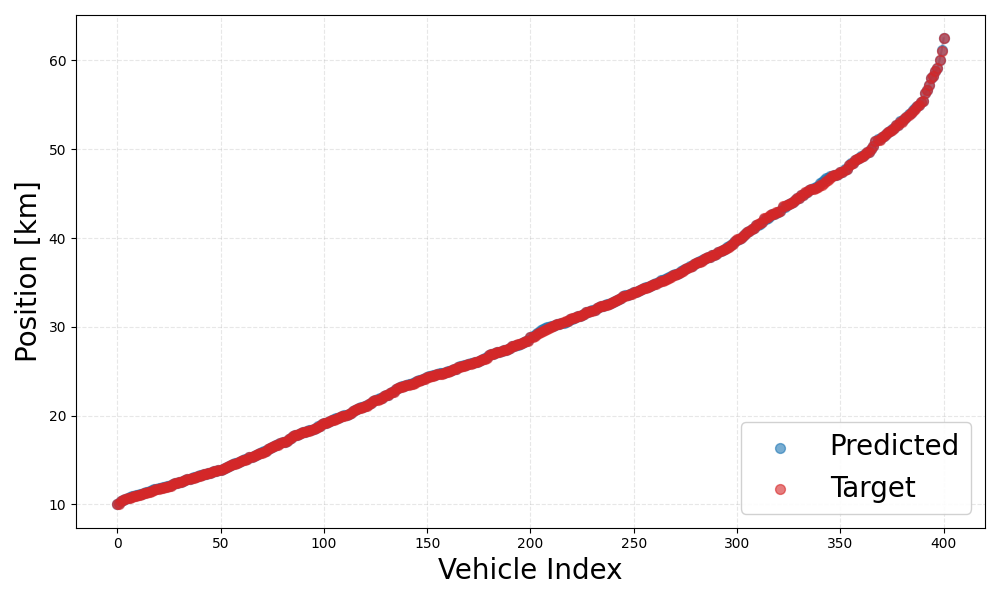}
	\caption{NN learning using 5,000 epochs}
	\label{fig:predicted_vs_target_34}
\end{subfigure}
\hfill
\begin{subfigure}[b]{0.23\textwidth}
	\centering
	\includegraphics[width=\textwidth]{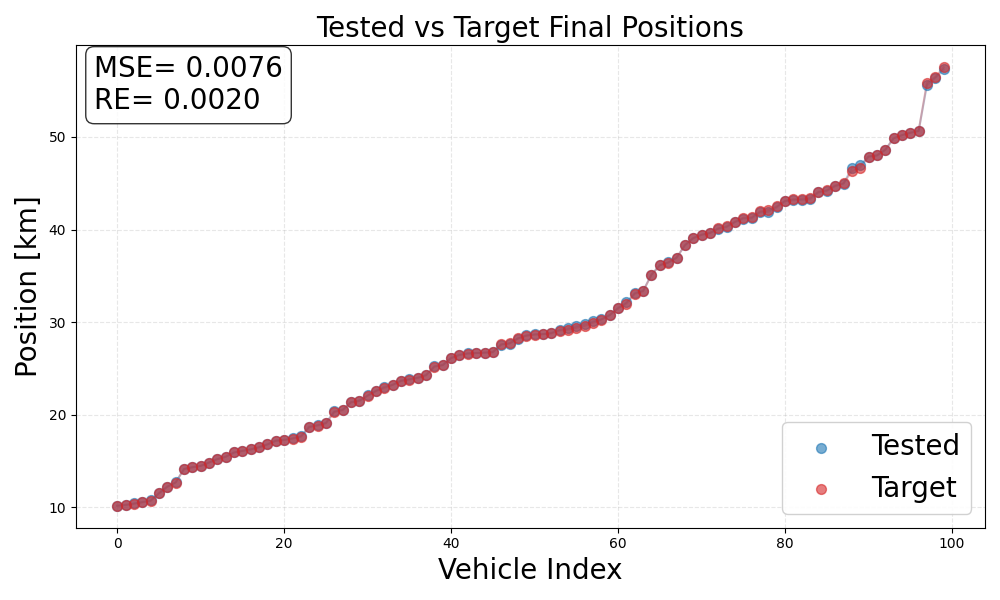}
	\caption{Testing using predicted density}
	\label{fig:tested_vs_target_34}
\end{subfigure}
\caption{Comparison between predictions and data with $N=4000$}
\vspace{-\baselineskip}
\label{fig:stop_and_go_recon_34}
\end{figure}

Fig.~\ref{fig:stop_and_go_recon_36} and Fig.~\ref{fig:stop_and_go_recon_34} illustrate the errors obtained by both training and testing phases of our learning procedure 
with different values of $N$ and $T$.
Precisely, during the testing phase where $N=2000$ and $T=0.1$, MSE was calculated to be $0.0483$, while RE was $0.0053$. When the number of vehicles was increased to $N=4000$ and the time duration to $T=0.2$, both MSE and RE decreased to $0.0076$ and $0.0020$ respectively. This demonstrates that our methods maintains its robustness even when applied to high dimensional data.

\begin{figure}[h]
\centering
\begin{subfigure}[b]{0.23\textwidth}
	\centering
	\includegraphics[width=\textwidth]{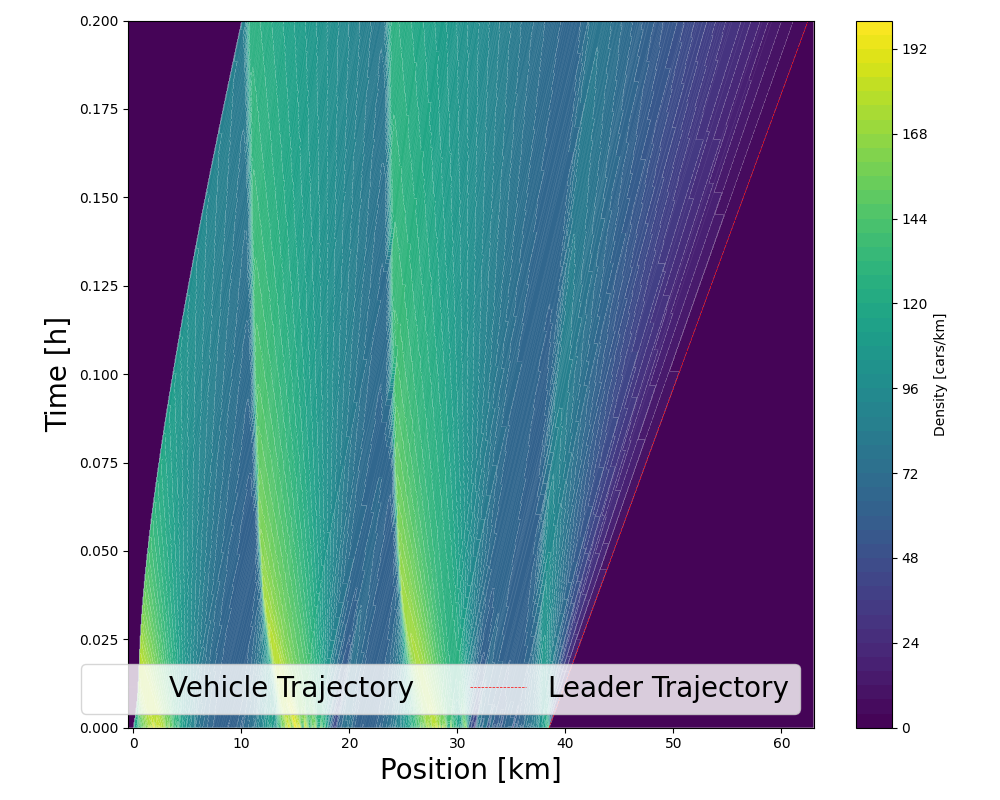}
	\caption{Reconstructed density}
\end{subfigure}
\hfill
\begin{subfigure}[b]{0.23\textwidth}
	\centering
	\includegraphics[width=\textwidth]{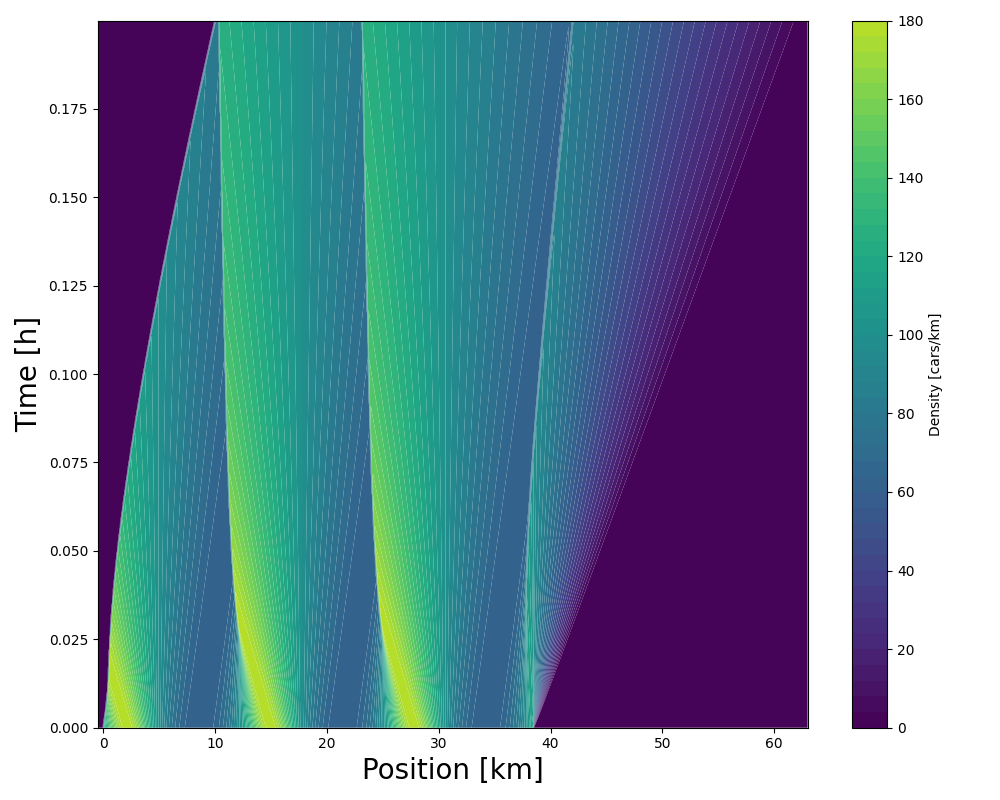}
	\caption{Godunov density}
\end{subfigure}
\caption{Comparison between the predicted density and Godunov scheme density with $N=4000$  and $T=0.2$}

\label{fig:stop_and_go_reconstruction_comparison}
\end{figure}

Additionally, the method allows us to reconstruct density for all times. Fig.~\ref{fig:stop_and_go_reconstruction_comparison} compares the reconstructed discrete density with the solution to LWR model \eqref{1}, which is computed using a Godunov scheme, as well as the ground truth initial density from which the total fleet of vehicles evolved before selecting probe vehicles. 


\begin{figure}
\begin{subfigure}[b]{0.23\textwidth}
	\centering
	\includegraphics[width=\textwidth]{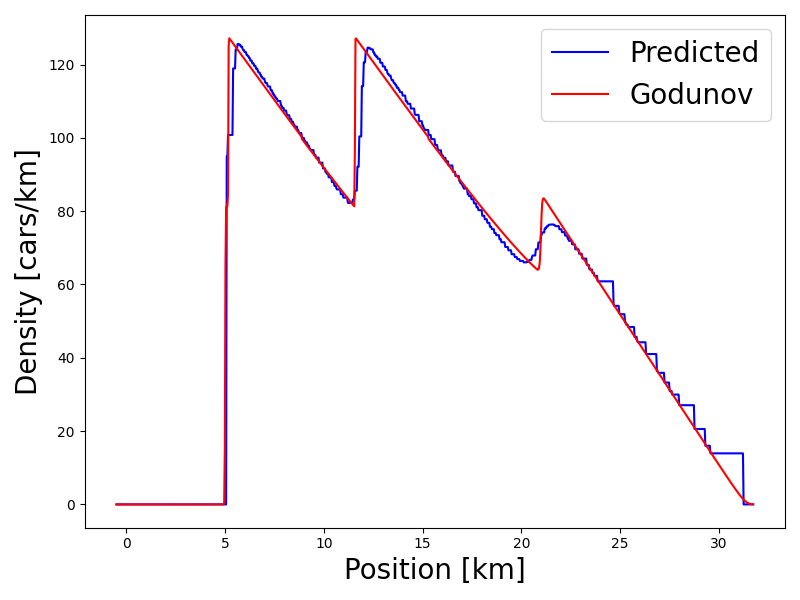}
	\caption{$N=2000$, $T=0.1$}
	\label{fig:final_densities_}
\end{subfigure}
\hfill
\begin{subfigure}[b]{0.23\textwidth}
	\centering
	\includegraphics[width=\textwidth]{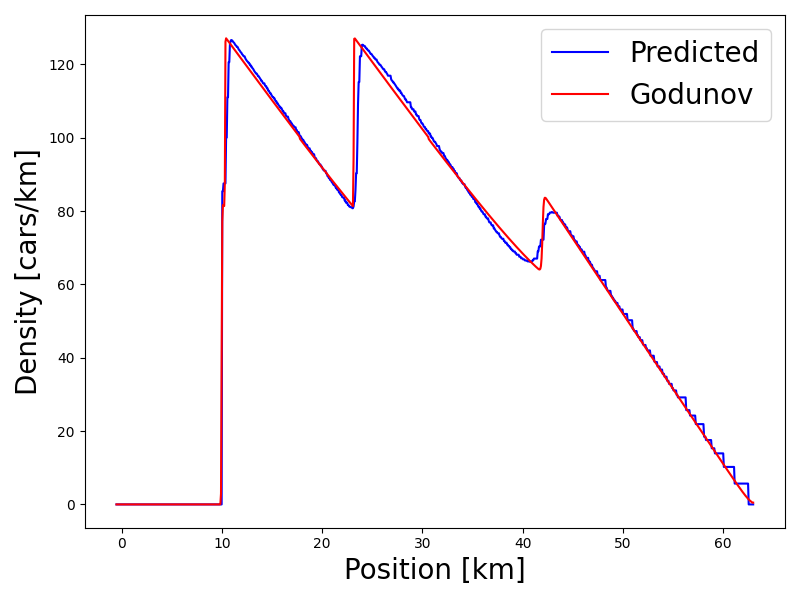}
	\caption{$N=4000$, $T=0.2$}
	\label{fig:final_densities_34}
\end{subfigure}
\caption{Reconstructed final traffic density}
\label{fig:final_stop_and_go_recon}
\end{figure}

Fig.~\ref{fig:final_stop_and_go_recon} presents the final reconstructed density for both previous settings, showing a strong match with the ground truth final density obtained by solving PDE \eqref{1}.

The second traffic scenario represents an abrupt transition in traffic conditions. The shock 
where occurs a change from low to high density is set at normalized position $0.5$ with respect to the initial domain, where position $0$ corresponds to the last vehicle and position $1$ corresponds to the leader's initial position. 
Let $\rho(0,x) = 
\begin{cases}
0.4\rho_{\max}, & 0 \leq x < 0.5, \\
0.9\rho_{\max}, & 0.5 < x \leq 1.
\end{cases}$
\begin{figure}
\begin{subfigure}[b]{0.23\textwidth}
	\centering
	\includegraphics[width=\textwidth]{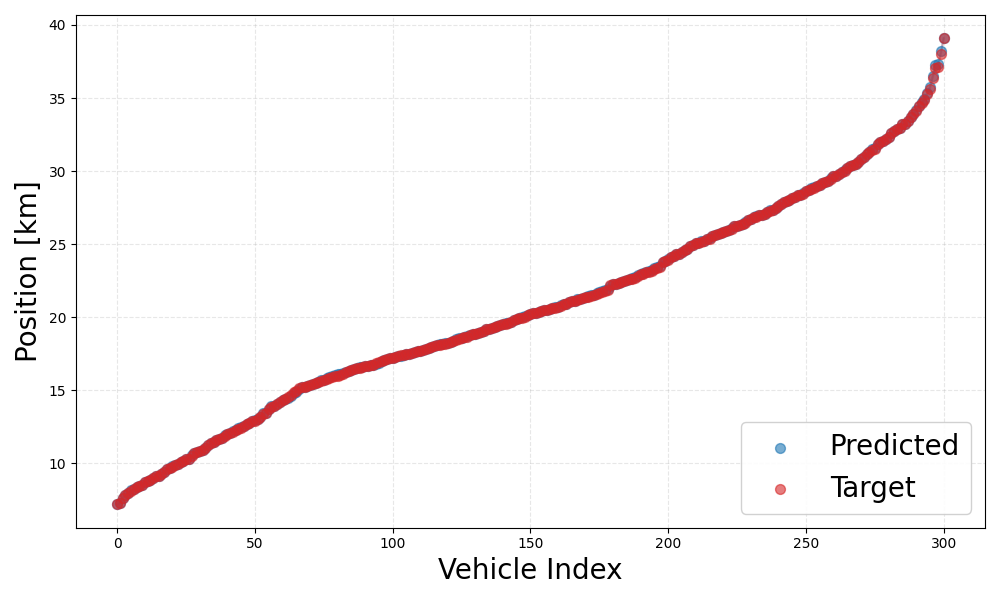}
	\caption{NN learning using 5,000 epochs}
	\label{fig:predicted_vs_target_30}
\end{subfigure}
\hfill
\begin{subfigure}[b]{0.23\textwidth}
	\centering
	\includegraphics[width=\textwidth]{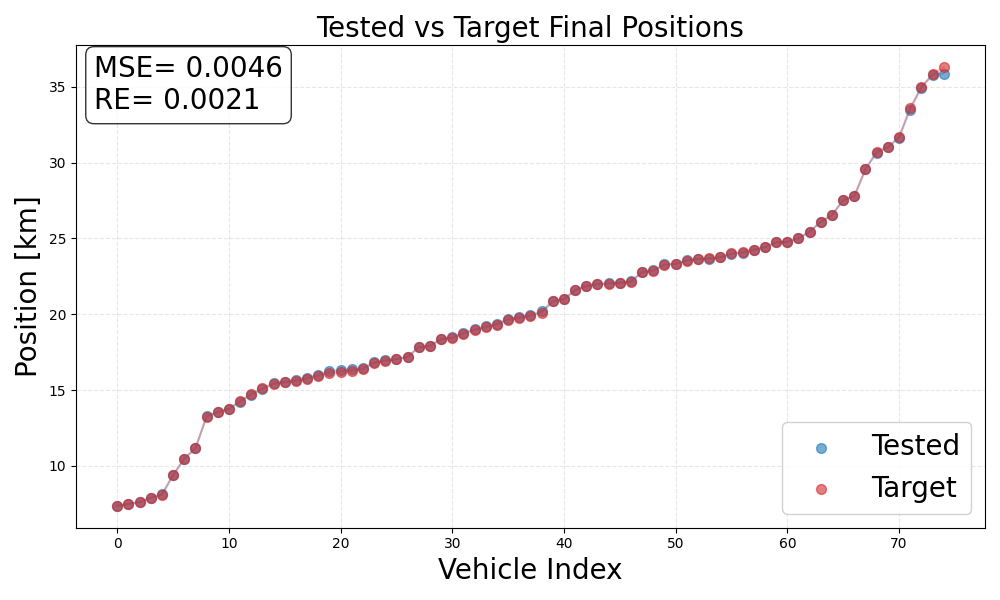}
	\caption{Testing using predicted density}
	\label{fig:tested_vs_target_30}
\end{subfigure}
\caption{Comparison between predictions and data with $N=3000$}
\vspace{-\baselineskip}
\label{fig:shock_recon_30}
\end{figure}



\begin{figure}[h]
\centering
\begin{subfigure}[b]{0.23\textwidth}
	\centering
	\includegraphics[width=\textwidth]{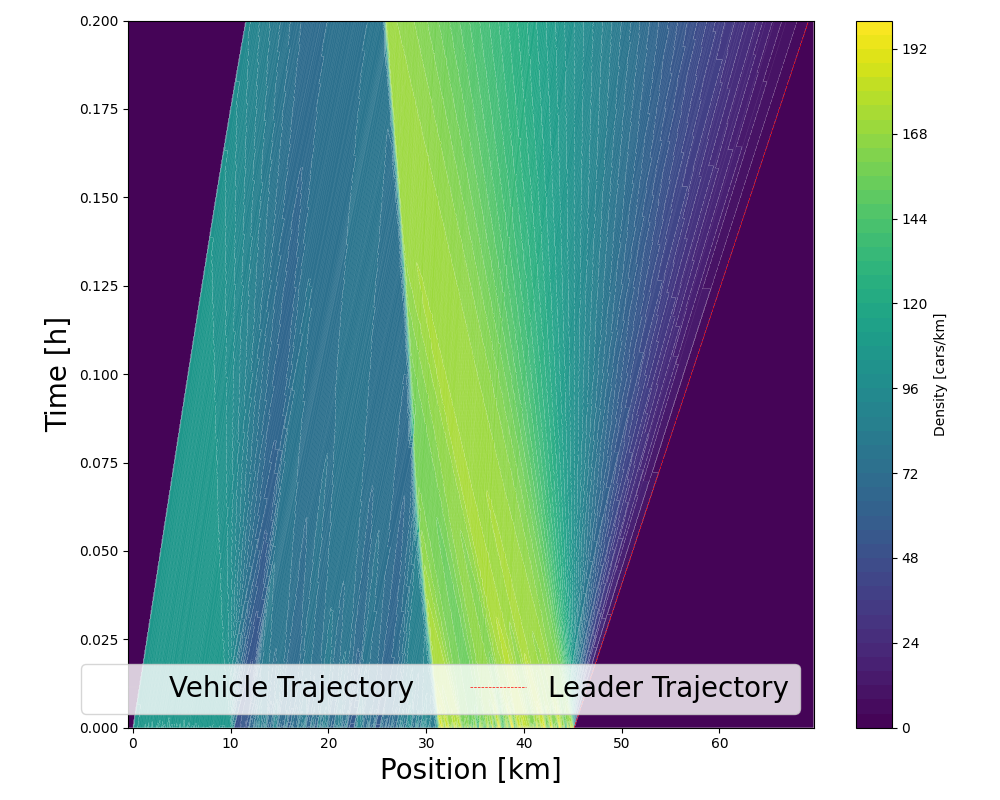}
	\caption{Reconstructed density}\label{fig:shock_reconstruction_comparison_reconstructed}
\end{subfigure}
\hfill
\begin{subfigure}[b]{0.23\textwidth}
	\centering
	\includegraphics[width=\textwidth]{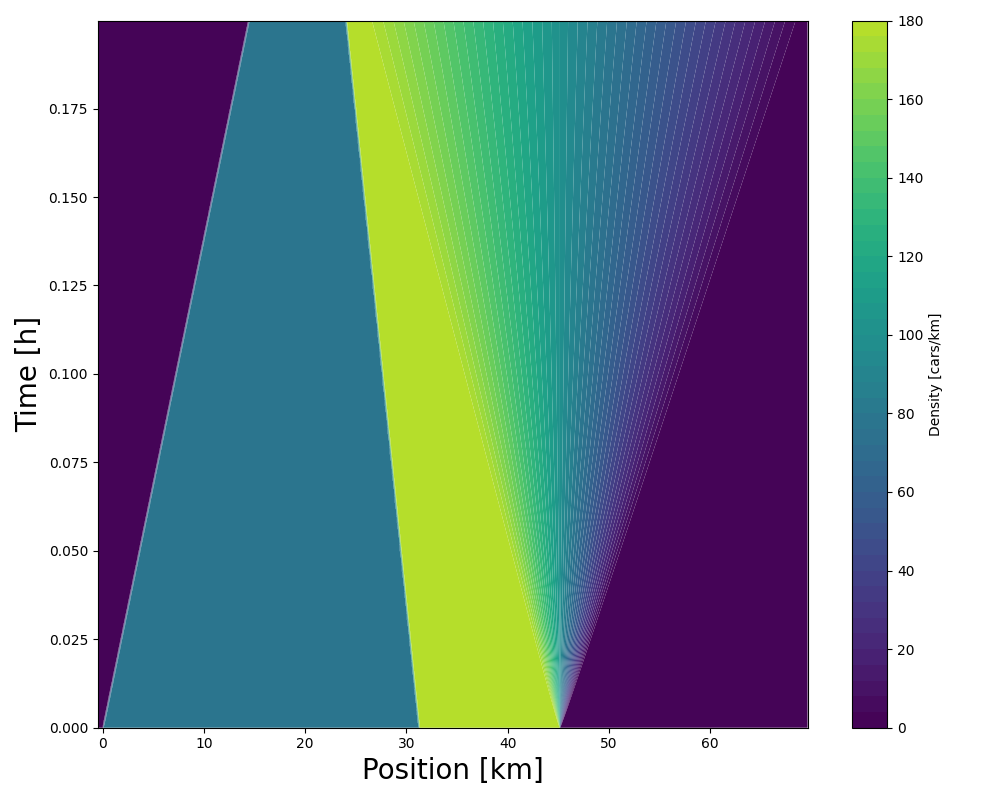}
	\caption{Godunov density}\label{fig:shock_reconstruction_comparison_godunov}
\end{subfigure}
\caption{Comparison between the predicted density and Godunov scheme density with $N=5000$  and $T=0.2$}
\label{fig:shock_reconstruction_comparison}
\end{figure}

\begin{figure}
\begin{subfigure}[b]{0.23\textwidth}
	\centering
	\includegraphics[width=\textwidth]{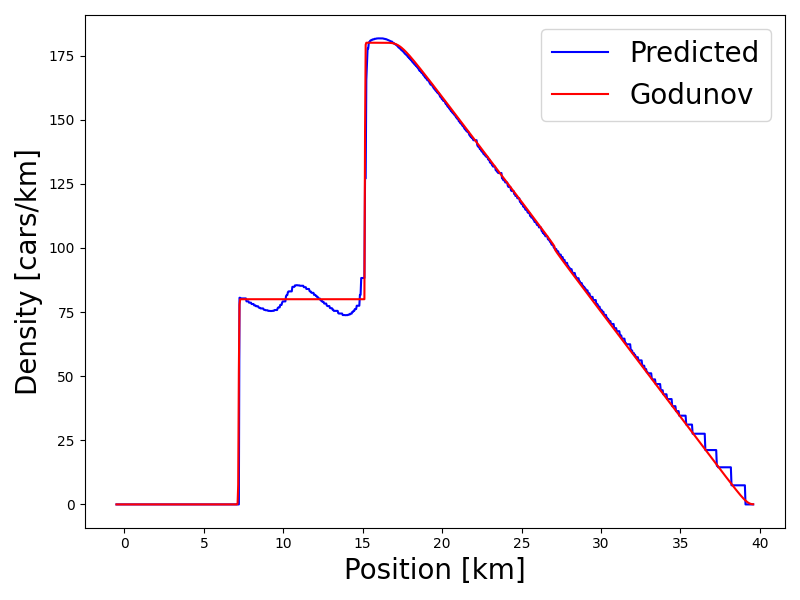}
	\caption{$N=3000$, $T=0.1$}
	\label{fig:final_densities}
\end{subfigure}
\begin{subfigure}[b]{0.23\textwidth}
	\centering
	\includegraphics[width=\textwidth]{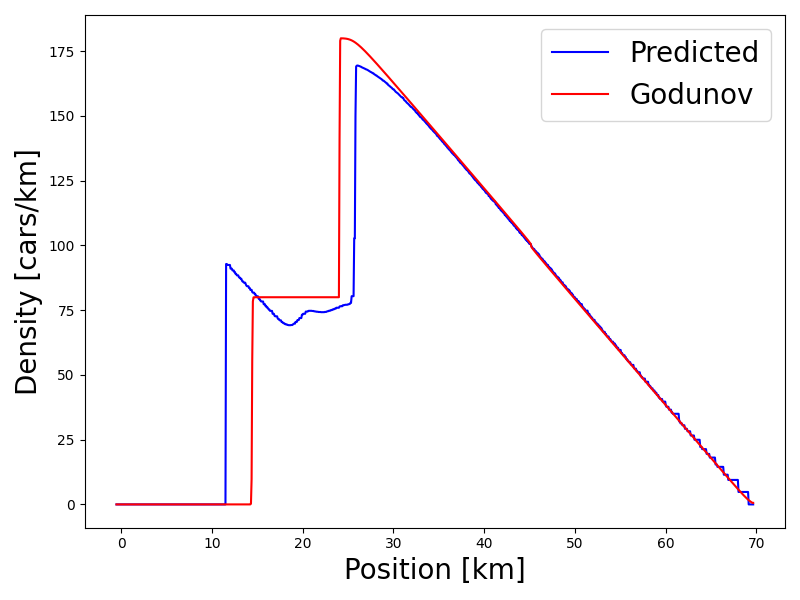}
	\caption{$N=5000$, $T=0.2$}
	\label{fig:final_densities_53}
\end{subfigure}
\caption{Reconstructed final traffic density}
\vspace{-\baselineskip}
\label{fig:final_shock_recon}
\end{figure}

The learning procedure applied to $N=3000$ cars yields results in Fig.~\ref{fig:shock_recon_30}.  In shock simulations, longer times reduce system variations, improving reconstruction consistency. Indeed, Fig.~\ref{fig:final_shock_recon} and Fig.~\ref{fig:shock_reconstruction_comparison} illustrate that for large final times and a large number of vehicles, such as $T=0.2$ and $N=5000$, the reconstructed final density converges well to the solution obtained from solving PDE \eqref{1}.  
Even under 
limited conditions, the reconstruction maintains an obvious accuracy. This suggests that the approach can effectively handle a range of scenarios, from those with ample time and vehicles to those with more limited resources. 


\section{Conclusions and Perspectives}
In this article, we have developed a model that combines traditional traffic flow models with a data-driven approach. Our method is able to perform effectively with limited data and sparse sensor information 
while maintaining low computational complexity. 
We avoid the need for real-time updates that often require high-frequency data inputs and can be expensive. We opted to use exclusively data generated by traffic FtL model \eqref{2b} rather than real-world data. The main reason lies on the fact that simulated data allow us to precisely replicate the assumptions and conditions required for theoretical guarantees. 
Our simulation focused on the use of artificial data generated by microscopic traffic models due to its suitability for theoretical convergence analysis. Indeed, we showed that simulated data allow us to precisely replicate the assumptions and conditions required for demonstrating the convergence of our model to the macroscopic LWR model. The methodology developed remains applicable to real-world data as a key strength of our model is its ability to integrate both artificial and real data. The use of real data would represent an important next step in validating our approach under realistic traffic conditions.	

\bibliographystyle{plain}
\bibliography{Biblio.bib}

@article{piccoli_vehicular_2009,
	title={Vehicular Traffic: A Review of Continuum Mathematical Models},
	author={Piccoli, B. and Tosin, A. and others},
	journal={Encyclopedia of Complexity and Systems Science},
	volume={22},
	pages={9727--9749},
	year={2009}
}

@article{bellomo_modeling_2011,
	title={On the Modeling of Traffic and Crowds: A Survey of Models, Speculations, and Perspectives},
	author={Bellomo, N. and Dogbe, C.},
	journal={SIAM Review},
	volume={53},
	number={3},
	pages={409--463},
	year={2011},
	publisher={SIAM}
}

@article{albi_vehicular_2019,
	title={Vehicular Traffic, Crowds, and Swarms: From Kinetic Theory and Multiscale Methods to Applications and Research Perspectives},
	author={Albi, G. and Bellomo, N. and Fermo, L. and Ha, S.-Y. and Kim, J. and Pareschi, L. and Poyato, D. and Soler, J.},
	journal={Math. Models Methods Appl. Sci.},
	volume={29},
	number={10},
	pages={1901--2005},
	year={2019},
	publisher={World Scientific}
}

@inproceedings{barreau_physics-informed_2021,
	author    = {Barreau, M. and Aguiar, M. and Liu, J. and Johansson, K. H.},
	title     = {Physics-Informed Learning for Identification and State Reconstruction of Traffic Density},
	booktitle = {Proceedings of the 60th IEEE Conference on Decision and Control (CDC)},
	year      = {2021},
	pages     = {2653--2658},
	publisher = {IEEE}
}

@article{colombo_well_2007,
	author={Colombo, R. M. and Goatin, P.},
	title={A Well-Posed Conservation Law with a Variable Unilateral Constraint},
	journal={J. Differential Equations},
	volume={234},
	number={2},
	pages={654--675},
	year={2007},
}

@article{francesco_rigorous_2017,
	author={Di Francesco, M. and Rosini, M. D.},
	title={Rigorous Derivation of Nonlinear Scalar Conservation Laws from Follow-the-Leader Type Models via Many Particle Limit},
	journal={Arch. Rational Mech. Anal.},
	volume={225},
	pages={129–166},
	year={2017}
}

@article{holden_continuum_2019,
	author={Holden, H. and Risebro, N. H.},
	title={The continuum limit of follow-the-leader models — A short proof},
	journal={SIAM J. Math. Anal.},
	volume={51},
	number={2},
	pages={1771--1786},
	year={2019}
}

@inproceedings{barreau_learning-based_2021,
	author    = {Barreau, M. and Liu, J. and Johansson, K. H.},
	title     = {Learning-Based State Reconstruction for a Scalar Hyperbolic {PDE} under Noisy Lagrangian Sensing},
	booktitle = {Proceedings of the 3rd Conference on Learning for Dynamics and Control},
	series    = {Proceedings of Machine Learning Research},
	volume    = {144},
	pages     = {34--46},
	publisher = {PMLR},
	year      = {2021}
}

@article{seo_traffic_2017,
	title={Traffic State Estimation on Highway: A Comprehensive Survey},
	author={Seo, T. and Bayen, A. M. and Kusakabe, T. and Asakura, Y.},
	journal={Annual Reviews in Control},
	volume={43},
	pages={128--151},
	year={2017},
	publisher={Elsevier}
}

@article{herrera_evaluation_2010,
	title={Evaluation of Traffic Data Obtained via GPS-Enabled Mobile Phones: The Mobile Century Field Experiment},
	author={Herrera, J. C. and Work, D. B. and Herring, R. and Ban, X. J. and Jacobson, Q. and Bayen, A. M.},
	journal={Transportation Research Part C: Emerging Technologies},
	volume={18},
	number={4},
	pages={568--583},
	year={2010},
	publisher={Elsevier}
}

@inproceedings{inzunza_pinn_2023,
	title={A PINN Approach for Traffic State Estimation and Model Calibration Based on Loop Detector Flow Data},
	author={Inzunza, D. and Goatin, P.},
	booktitle={8th International Conference on Models and Technologies for Intelligent Transportation Systems (MT-ITS)},
	year={2023},
	pages={1--6}
}

@article{shi_physics_2021,
	author = {Shi, R. and Mo, Z. and Huang, K.  and Di, X. and Du, Q.},
	title = {A Physics-Informed Deep Learning Paradigm for Traffic State and Fundamental Diagram Estimation},
	journal = {IEEE Transactions on Intelligent Transportation Systems},
	year = {2021},
}
\end{document}